\documentclass[12pt]{amsart}

\usepackage{amssymb}
\usepackage{amsmath}
\usepackage{amsthm}

\input xy
\xyoption{all}

\newcommand{\ran}{\mathop{\mathrm{ran}}}

\title{Finite dimensional point derivations for graph algebras}

\author{Benton L. Duncan}

\address{Department of Mathematics\\
300 Minard Hall\\
North Dakota State University\\
Fargo, ND  58105-5075\\
USA}

\email{benton.duncan@ndsu.edu}

\subjclass[2000]{47L40, 47L55, 47L75, 46L80}

\keywords{noncommutative point derivations, graph operator algebras}

\begin{document}

\theoremstyle{plain}
\newtheorem{thm}{Theorem}[section]
\newtheorem{lem}[thm]{Lemma}
\newtheorem{prop}[thm]{Proposition}
\newtheorem{cor}[thm]{Corollary}

\theoremstyle{definition}
\newtheorem{dfn}[thm]{Definition}
\newtheorem*{construction}{Construction}
\newtheorem*{example}{Example}

\theoremstyle{remark}
\newtheorem*{question}{Question}
\newtheorem*{rmk}{Remark}

\begin{abstract}  This paper focuses on certain finite dimensional point
derivations for the non-selfadjoint operator algebras corresponding
to directed graphs. We begin by analyzing the derivations
corresponding to full matrix representations of the tensor algebra
of a directed graph.  We determine when such a derivation is inner,
and describe situations that give rise to non-inner derivations. We
also analyze the situation when the derivation corresponds to a
multiplicative linear functional.
\end{abstract}

\maketitle

\section{Introduction}

The non-selfadjoint operator algebras associated to directed graphs
have undergone significant scrutiny of late.  Specifically, the
algebra $\mathcal{T}^{+}(Q)$, the norm closed algebra generated by
the left regular representation of a directed graph $Q$, and
$\mathcal{L}_Q$, the WOT-closure of $\mathcal{T}^{+}(Q)$, are
studied abstractly as special cases of tensor algebras over
$C^*$-correspondences in \cite{Muhly-Solel:1998}. The case of tensor
algebras corresponding to directed graphs were first studied in
\cite{Pop:1996} in the case of the graph with a single vertex and
$n$-edges.  The weakly closed version was studied around the same
time in \cite{Dav-Pitts:1998, Dav-Pitts:1998a}. General directed
graphs were taken up in \cite{Kribs-Power:2003a}. Since that time
the study has expanded significantly to various facets of these
non-selfadjoint algebras, see for example
\cite{Davidson-Katsoulis:2004, Jaeck-Power:2005, Jury-Kribs:2005,
Kat-Kribs:2003, Kribs-Power:2003b}.

Important in the analysis of these algebras is a recent paper
\cite{Davidson-Katsoulis:2004} which analyzed finite dimensional
representations and faithful irreducible representations for
strongly transitive directed graphs. Building on their work, and
combining with a recent paper \cite{Duncan:2005} we undertake a
description of noncommutative point derivations of directed graph
algebras.  When one views the classical situation of uniform
algebras a strong connection is found between the point derivations
and analytic structure \cite{Browder:1986}. In particular for the
disk algebra there exist point derivations at a character if and
only if the character corresponds to a point on the interior of the
unit disk, see Section 9 of \cite{Johnson:1972}. In the present
paper we find similar results when one looks at point derivations at
characters.

On the other hand when dealing with noncommutative algebras it is
clear that the characters are not enough to say much about the
algebra.  We have thus taken a look at point derivations at certain
irreducible finite dimensional representations of directed graph
operator algebras.  Here too we find a sense of analytic structure
but we have not developed that theory to its full extent.

This work is developed as further strengthening of the connection
between the disk algebra and directed graph operator algebras which
was begun in \cite{Duncan:2005} in the case of a directed cycle
graph. There we developed the theory in direct analogy with the
standard results for uniform algebras.  Here, however we take a
different approach, since the finite dimensional representations are
more diverse. More importantly, the general situation of
$\mathcal{T}^{+}(Q)$ and $\mathcal{L}_Q$-valued derivations of
$\mathcal{T}^{+}(Q)$ are not amenable to the approach of
\cite{Duncan:2005} where the finite dimensional representations are
used to make statements about general derivations.

In the first part of the paper we develop a description of a rich
class of finite dimensional representations, generalizing a specific
case from \cite{Davidson-Katsoulis:2004}.  We then factor our
representations through the algebras given by the directed cycle
graphs. It is then a simple application of the results of
\cite{Duncan:2005} to describe the noncommutative point derivations
into $M_n$ for general directed graph algebras.  This suggests a
notion of noncommutative analyticity coming from the finite
dimensional representations.

We take a similar approach to studying the point derivations at
characters for a general directed graph operator algebra.  Once
again we factor such a representation through the a well understood
example, the noncommutative analytic Toeplitz algebras of
\cite{Dav-Pitts:1998}.  Thus the discussion of point derivations and
commutative analytic structure reduces to the discussion of point
derivations on these well understood graph algebras.

We close the paper with a result concerning the range of a
derivation $D: A_n \rightarrow A_n$.  We show that such a derivation
must have range contained in the commutator ideal of $A_n$. Of
course an inner derivation will satisfy this property. However our
result says nothing to suggest that every derivation on $A_n$ is
inner.  In fact we have little evidence at this point that such a
result is even true.

\section{Notation and Background}

We begin with a review of background material and we fix some
notation.  To a directed graph $Q$ there exists two non-selfadjoint
operator algebras which we will study below. Both arise from the
left regular representation of the graph acting on $\ell^2$ of the
finite path space.  The first algebra, $\mathcal{T}^{+}(Q)$ will be
the norm closure of this representation.  The second $\mathcal{L}_Q$
will be the WOT closure of this representation.

For a directed graph $Q$ we denote the edge set of $Q$ by $E(Q)$ and
the vertex set of $Q$ by $V(Q)$.  To each edge there are maps $r:
E(Q) \rightarrow V(Q)$ and $ s: E(Q) \rightarrow V(Q)$ which give
the range and source of an edge, respectively.  We will write
$\mathcal{C}_n$ for the cycle graph given by $n$ distinct edges $\{
e_i \}$ and $n$ vertices $\{ v_i \}$, with $s(e_i) = r(e_{i+1})$ for
$ 1 \leq i \leq n-1$ and $ s(e_n) = r(e_1)$.  Recall that this
algebra can be written as a matrix function algebra of the form \[
\begin{bmatrix} f_{1,1}(z^n) & z f_{1,2}(z^n) & z^2
f_{1,3}(z^n) & \cdots & z^{n-1}f_{1,n}(z^n) \\ z^{n-1} f_{2,1}(z^n)
& f_{2,2}(z^n) & z f_{2,3}(z^n) & \cdots & z^{n-2}f_{2,n}(z^n) \\
z^{n-2}f_{3,1}(z^n) & z^{n-1}f_{3,2}(z^n) & f_{3,3}(z^n) & \cdots &
z^{n-3}f_{3,n}(z^n) \\ \vdots & \vdots & \vdots & \ddots & \vdots
\\ zf_{n,1}(z^n) & z^2 f_{n,2}(z^n) & z^{3}f_{n,3}(z^n) &
\cdots & f_{n,n}(z^n) \end{bmatrix} \] where $f_{i,j} \in
A(\mathbb{D})$ for all $ 1\leq i,j \leq n$.

We will denote by $B_n$ the graph with $1$ vertex and $n$ edges. For
shorthand the algebras, $\mathcal{T}^{+}(B_n)$ and
$\mathcal{L}_{B_n}$ will be denoted by $A_n$ and $\mathcal{L}_n$,
respectively.  In this case we will assign an ordering to the edges
and denote the isometry associated to the $i$-th edge by $L_i$.

We now establish some standard definitions and notation.  A path in
$Q$ will be finite sequence $e_1e_2, \cdots e_n$ with $e_i \in E(Q)$
and $ r(e_i) = s(e_{i-1})$ for all $ 2 \leq i \leq n$. Recall that a
directed graph is {\em transitive} if for every pair of vertices,
$v$ and $w$ there is a directed path beginning at $v$ and ending at
$w$. We say that a path $w$ in a directed graph is {\em primitive}
if $w \neq v^n$ for any paths $v$.  We say that a path $w = e_1 e_2
\cdots e_n$ is a {\em cycle} if $ r(e_1) = s(e_n)$.

We will write $M_n$ to mean the $n \times n$ matrices with entries
from $\mathbb{C}$.  We denote by $e_{i,j}$ the elementary matrix in
$M_n$ with $1$ in the $i$-$j$ position and $0$ everywhere else.

We now state a standard result concerning representations of graph
operator algebras which will be useful in what follows.  This result
follows from work in \cite{Muhly-Solel:1998}, or
\cite{Jury-Kribs:2004} and is given explicitly in the case of
countable directed graphs in Section 3 of
\cite{Davidson-Katsoulis:2004}.

\begin{prop}\label{ccrepn} Let $Q$ be a directed graph and $A$ an
operator algebra.  Let $\pi: Q \rightarrow A$ be a map with $\pi(v)
=P_v$ a projection for all $v \in V(Q)$, $\pi(e)=L_e $ is a nonzero
contraction for all edges $e \in E(Q)$, and:
\begin{enumerate} \item $P_v$ is orthogonal to $P_w$ for all
vertices $v, w \in V(Q)$.

\item $P_{r(e)}L_eP_{s(e)} = L_e$ for all edges $e \in E(Q)$.

\item $[L_e]_{e \in E(Q)}$ is a row contraction in $A$.
\end{enumerate}
Then $\pi$ extends to a completely contractive representation of
$\mathcal{T}^{+}(Q)$. \end{prop}

\section{Irreducible representations into $M_n$}

Let $Q$ be a directed graph and let $ w = e_1e_2 \cdots e_n$ be a
finite path in $Q$.  For $\lambda \in \overline{\mathbb{D}}$, $\mu
\in \mathbb{T}$, $v \in V(Q)$, and $e \in E(Q)$ define \[ \pi_{w,
\lambda, \mu}(P_v) = \sum_{s(e_{j}) = v} e_{j,j} \mbox{ and }
\pi_{w, \lambda, \mu}(L_e) = \sum_{e_{j} = e} \lambda e_{j-1,j} \]
where for the sake of the notation we denote by $e_{0,1}$ the matrix
$\mu e_{n,1}$. This map then extends to a representation $\pi_{w,
\lambda, \mu}: \mathbb{F}^{+}(Q) \rightarrow M_n$.

\begin{prop}\label{primitive} Let $ w$ be a finite path in $Q$, $\lambda
\in \overline{\mathbb{D}}$, and $\mu \in \mathbb{T}$, then the
representation $\pi_{w, \lambda, \mu}: \mathbb{F}^{+}(Q) \rightarrow
M_n$ extends to a completely contractive representation of
$\mathcal{T}^{+}(Q)$ into $M_n$.  Moreover, if $w$ is a primitive
cycle, and $ \lambda \neq 0$, then the extension is onto.\end{prop}

\begin{proof} The extension to a completely contractive
representation of the algebra $\mathcal{T}^{+}(Q)$ follows from
Proposition \ref{ccrepn}.  We discuss the other conclusion now.  The
details of the argument are in the proof of \cite[Lemma
4.3]{Davidson-Katsoulis:2004}.  We only sketch the proof here.

Let $w = e_1e_2\cdots e_n$ and notice that \[\pi_{w, \lambda,
\mu}(L_{e_je_{j-1}e_{j-2} \cdots e_1 e_n e_{n-1} \cdots e_j}) =
\lambda^{-k-1} \mu e_{j,j+1}.\]  Since $ \lambda \neq 0$, letting
$j$ vary yields a generating set for $M_n$ in the range of
$\mathcal{T}^{+}(Q)$ and hence the extension is onto. \end{proof}

\begin{rmk} If $ |\lambda|<1$ then the map is $w^*$-continuous by
Corollary 3.2 in \cite{Davidson-Katsoulis:2004} so that the
representation extends to a $w^*$-continuous completely contractive
representation of $\mathcal{L}_Q$. \end{rmk}

Notice that in the special case of $\lambda = \frac{1}{2}$ this
representation is the representation $\varphi_{w, \mu}$ given in
Section 4 of \cite{Davidson-Katsoulis:2004}, with respect to the
usual orthonormal basis of $\mathbb{C}^n$ with reverse ordering.  It
follows, by \cite[Theorem 4.4]{Davidson-Katsoulis:2004} that in the
case of transitive graphs, by letting $w, \varphi,$ and $\mu$ vary,
we get a family of irreducible representations which separate the
points of $\mathcal{T}^{+}(Q)$.  Here however the finite dimensional
representations are richer and will allow for a more detailed
discussion of the noncommutative point derivations.

\begin{dfn}  Let $A$, $B$, and $C$ be operator algebras. Let
$\pi: A \rightarrow B$ be a completely contractive representation of
the operator algebra $A$.  We say that $\pi$ {\em factors through
$C$} if there exist completely contractive representations, $\iota:
A \rightarrow \mathbb{C}$ and $\tilde{\pi}: C \rightarrow B$ such
that $ \pi(a) = \tilde{ \pi} \circ \iota(a)$ for all $a \in A$.
\end{dfn}

We will show that the representations $\pi_{w, \lambda, \mu}:
\mathcal{T}^{+}(Q) \rightarrow M_n$ factor through
$A(\mathbb{C}_n)$.  We begin by constructing the map $\iota:
\mathcal{T}^{+}(Q) \rightarrow A(\mathbb{C}_n)$. For notation sake,
with $1 \leq i <n $ let $Z_i \in \mathcal{T}^{+}(\mathcal{C}_n)$ be
the matrix with $z$ in the $i$-$(i+1)$ position and zeroes
everywhere else.  Denote by $Z_n$ the matrix in
$\mathcal{T}^{+}(\mathcal{C}_n)$ with a $z$ in the $n$-$1$ position
and zeroes everywhere else.

For a finite path $w$ in $Q$ given by $w = e_1e_2 \cdots e_n$ define
the representation $\iota_w: \mathbb{F}^{+}(Q) \rightarrow
\mathcal{T}^{+}(\mathcal{C}_n)$ by first setting, for $v \in V(Q)$,
\[ \iota_w(P_v) = \sum_{s(e_j) = v} e_{j,j} \in \mathcal{T}^{+}(\mathcal{C}_n). \]
Next, for $e \in E(Q)$ define \[ \iota_w(L_e) = \sum_{e_j = e}
Z_{j}. \] The map $\iota_w$ will then be the natural extension to
$\mathbb{F}^{+}(Q)$.  The next proposition follows immediately from
Proposition \ref{ccrepn}.

\begin{prop} Given a finite path $w$, the map $\iota_w$ extends to a
completely contractive representation $\iota_w: \mathcal{T}^{+}(Q)
\rightarrow \mathcal{T}^{+}(\mathcal{C}_n)$. \end{prop}

\begin{rmk} This map is also $w^*$-continuous and sends
$\mathcal{L}_Q$ into $\mathcal{L}_{\mathcal{C}_n}$.  This also
follows from Corollary 3.2 of \cite{Davidson-Katsoulis:2004}, by
noticing that the left regular representation of
$\mathcal{L}_{\mathcal{C}_n}$ is pure. \end{rmk}

An easy consequence of the definition is the following lemma.

\begin{lem} For a primitive cycle $w$, the map $\iota_w: \mathcal{T}^{+}(Q)
\rightarrow \mathcal{T}^{+}(\mathcal{C}_n)$ is onto if and only if
$s(e_i) \neq s(e_j)$ for all $ i \neq j$.\end{lem}

\begin{proof} If $\iota_w$ is onto, then for each $i$ there is
$X_{i,i} \in \mathcal{T}^{+}(Q)$ with $\iota_w(X_{i,i}) = e_{i,i}$.
But notice that by definition $ \iota_w(X_{i,i}) = e_{i,i}$ if and
only if $X_{i,i} = P_{v_i}$ where $ v_i = s(e_i)$ and $ v_i \neq
v_j$ for all $i \neq j$.

Notice that if $s(e_i) \neq s(e_j)$ for all $i \neq j$ then in
particular we know that $e_i \neq e_j$ for all $i \neq j$. It
follows that $\iota_w(L_{e_i}) = Z_i$.  Similarly
$\iota_w(P_{s(e_i)}) = e_{i+1,i+1}$ and hence the range of $\iota_w$
contains a generating set for
$\mathcal{T}^{+}(\mathcal{C}_n)$.\end{proof}

Notice that even when $\iota_w$ is not onto the same argument as in
the proof of Proposition \ref{primitive} will tell us that if
$A_0(z^n)$ is the nonunital subalgebra of $A(z^n)$ generated by
$z^n$ then \[ \begin{bmatrix} A_0(z^n)& z A_0(z^n) & z^2 A_0(z^n) &
\cdots & z^{n-1}A_0(z^n) \\ z^{n-1} A_0(z^n)
& A_0(z^n) & z A_0(z^n) & \cdots & z^{n-2}A_0(z^n) \\
z^{n-2}A_0(z^n) & z^{n-1}A_0(z^n) & A_0(z^n) & \cdots &
z^{n-3}A_0(z^n) \\ \vdots & \vdots & \vdots & \ddots & \vdots
\\ zA_0(z^n) & z^2 A_0(z^n) & z^{3}A_0(z^n) &
\cdots & A_0(z^n) \end{bmatrix}  \subseteq \ran (\iota_w). \] It
follows that the complement of $\ran \iota_w$ is a finite
dimensional subspace of $\mathcal{T}^{+}(\mathcal{C}_n)$ and hence
is complemented.  This will be important when we construct
derivations.

Now, for $ \lambda \in \overline{\mathbb{D}}$ and $ \mu = e^{i
\theta} \in \mathbb{T}$ we define a completely contractive
representation $\pi_{\lambda, \mu}: \mathcal{T}^{+}(\mathcal{C}_n)
\rightarrow M_n$. The map $ \pi_{\lambda, \mu}$ will be chosen so
that $\pi_{\lambda, \mu}(e_{i,i}) = e_{i,i}$, $\pi(Z_i) = \lambda
e_{i,i+1}$ for $ 1 \leq i < n$, and $\pi(Z_n) = \mu \lambda
e_{n,1}$.  We begin this by noticing that the inner automorphism
$\pi_{\mu}: A(\mathbb{C}_n) \rightarrow A(\mathbb{C}_n)$ induced by
the matrix \[
\begin{bmatrix} e^{i \frac{\theta}{n}} & 0 & 0 & \cdots & 0 & 0 \\ 0
& e^{-i \frac{\theta(n-2)}{n}} & 0 * \cdots & 0 & 0 \\ 0 & 0 & e^{-i
\frac{\theta(n-3)}{n}} & \cdots & 0 & 0 \\ \vdots & \vdots & \vdots
& \ddots & \vdots & \vdots \\ 0 & 0 & 0 & \cdots & e^{-i
\frac{\theta}{n}} & 0 \\ 0 & 0 & 0 & \cdots & 0 & 1 \end{bmatrix}.
\]  A technical calculation tells us that for all $i$, $ \pi_{\mu}
(e_{i,i}) = e_{i,i}$, $ \pi_{\mu}(Z_i) = Z_i$ for $ 1 \leq i < n$
and $ \pi_{\mu}(Z_n) = \mu Z_n$.  Following the automorphism by the
completely contractive representation $\tau_{\lambda}:
\mathcal{T}^{+}(\mathcal{C}_n) \rightarrow M_n$, which is given by
evaluation at $\lambda$, yields a map $ \pi_{\lambda, \mu}$ as
described.

\begin{thm} For $ \lambda \in \overline{\mathbb{D}}$, $\mu
\in \mathbb{T}$ and $w$ a finite path in $Q$, the representation $
\pi_{w, \lambda, \mu} : \mathcal{T}^{+}(Q) \rightarrow M_n$ factors
through $ \mathcal{T}^{+}(\mathcal{C}_n)$ via the representation $
\tau_{\lambda} \circ \pi_{\mu} \circ \iota_w$.\end{thm}

\begin{proof} Notice that $\iota_w$ is completely contractive by
construction.  Further, $\pi_{\mu}$ is completely contractive as it
is given as an inner automorphism by an invertible of norm $1$. For
$\tau_{\lambda}$ notice that the automorphism $a_{\lambda}:
A(\mathbb{D}) \rightarrow A(\mathbb{D})$ given by $a_{\lambda}(f(z))
= f(\lambda)$ is contractive, and hence completely contractive.
Thus, the matricial version of the automorphism $a_{\lambda}^{(n)}:
M_n \otimes A(\mathbb{D}) \rightarrow M_n \otimes A(\mathbb{D})$ is
completely contractive.  Notice however, that $\tau_{\lambda}$ is
the restriction to $\mathcal{T}^{+}(\mathcal{C}_n)$ of the map
$a_{\lambda}^{(n)}$, and hence $\tau_{\lambda}$ is completely
contractive.

It follows that the composition, $\tau_{\lambda} \circ \pi_{\mu}
\circ \iota_w$ is a completely contractive representation of
$\mathcal{T}^{+}(Q)$ into $M_n$.  Notice that $\tau_{\lambda} \circ
\pi_{\mu} \circ \iota_w (L_e) = \pi_{w, \lambda, \mu} (L_e)$ for all
$e \in E(Q)$ and $\tau_{\lambda} \circ \pi_{\mu} \circ \iota_w (P_v)
= \pi_{w, \lambda, \mu} (P_v)$ for all $v \in V(Q)$.  It follows
that, as $\mathcal{T}^{+}(Q)$ is generated by \[ \{P_v, L_e: v \in
V(Q), e \in E(Q) \},\] we know that $\tau_{\lambda} \circ \pi_{\mu}
\circ \iota_w (a) = \pi_{w, \lambda, \mu} (a)$ for all $ a \in
\mathcal{T}^{+}(Q)$.
\end{proof}

\begin{rmk} Notice that $\pi_{\mu}$ and $ \iota_w$ are completely
contractive and $w^*$-continuous.  Thus if $|\lambda|<1$ we know
that the $w^*$-continuous extension of $\pi_{w, \lambda, \mu}$ to
all of $\mathcal{L}_{Q}$ factors through
$\mathcal{L}_{\mathcal{C}_n}$ via $w^*$-continuous representations.
\end{rmk}

\section{Noncommutative point derivations into $M_n$}

Given a completely contractive representation $\pi: A \rightarrow B$
we say that a continuous linear map $D: A \rightarrow B$ is a {\em
derivation at $\pi$} if $ D(ab) = D(a) \pi(b) + \pi(a) D(b)$ for all
$ a, b \in A$.  It was shown in \cite[Proposition 1]{Duncan:2005}
that if $x \in B$ the linear map $\delta_X(a) = \pi(a) X- X \pi(a)$
for all $ a \in A$ is a derivation at $\pi$. Any derivation at $\pi$
of this form is said to be {\em inner at $\pi$}.

For $w$ a finite primitive cycle, $\lambda \in
\overline{\mathbb{D}}$, and $ \mu \in \mathbb{T}$ we will be
interested in continuous linear maps $D: \mathcal{T}^{+}(Q)
\rightarrow M_n$ which are derivations at $ \pi_{w,\lambda,\mu}$. We
first find a method to recognize whether a derivation is inner.

\begin{lem} Let $Q$ be a directed graph, $w$ a primitive cycle
in $Q$, $\lambda \in \overline{\mathbb{D}}$ with $ \lambda \neq 0$,
and $ \mu \in \mathbb{T}$. For $D: \mathcal{T}^{+}(Q) \rightarrow
M_n$, a continuous derivation at $ \pi_{w, \lambda, \mu}$, $D$ is
inner if and only if $D(a) = 0 $ for all $ a \in \ker( \pi_{w,
\lambda, \mu})$. \end{lem}

\begin{proof}  If $D$ is inner there exists, by definition, an $X \in
M_n$ such that $D(a) = \pi_{w, \lambda, \mu}(a) X - X \pi_{w,
\lambda, \mu}(a)$ for all $ a \in \mathcal{T}^{+}(Q)$. Now if $a \in
\ker \pi_{w, \lambda, \mu}$ then $D(a) = 0X-X0 = 0$ and hence
$D|_{\ker( \pi_{w, \lambda, \mu})} = 0$.

Now assume that $D(a) = 0$ for all $a \in \ker ( \pi_{w, \lambda,
\mu})$. Define a map $\widehat{D}: M_n \rightarrow M_n$ by
$\widehat{D}(x) = D(a)$ where $ \pi_{w, \lambda, \mu}(a) = x$.
Notice that if $\pi_{w, \lambda, \mu}(a) = \pi_{w, \lambda, \mu}(b)
= x$ then $ a-b \in \ker \pi_{w, \lambda, \mu}$ and hence $D(a-b) =
0$.  Thus $\widehat{D}$ is well defined.  Since $\pi_{w, \lambda,
\mu}$ is onto, for every $x,y \in M_n$ we have $a, b \in
\mathcal{T}^{+}(Q)$ such that $\pi_{w, \lambda, \mu}(a) = x$ and $
\pi_{w, \lambda, \mu}(b) = y$. Notice that $\widehat{D}(xy) = D(ab)$
by definition, but $D(ab) = D(a) \pi_{w, \lambda, \mu}(b) + \pi_{w,
\lambda, \mu}(a) D(b) = \widehat{D}(x)y + x \widehat{D}(y)$ and
hence $\widehat{D}$ is a derivation on $M_n$.

It is well known that every derivation on $M_n$ is inner and hence
there is $X \in M_n$ with $ \widehat{D}(y) = yX-Xy$ for all $x \in
M_n$. By definition, for any $a$ in $\mathcal{T}^{+}(Q)$ with $
\pi_{w, \lambda,
\mu}(a) = y$ we know that \begin{align*} D(a) &= \widehat{D}(y) \\
&= yX-Xy \\ &= \pi_{w, \lambda, \mu}(a) X - X \pi_{w, \lambda,
\mu}(a) \end{align*} and hence $D$ is inner.\end{proof}

We now look at derivations in the special case where $Q$ is the
graph $\mathcal{C}_n$, as developed in \cite{Duncan:2005}.

\begin{prop} Let $w$ be a primitive cycle in $\mathcal{C}_n$ and
$\lambda \in \overline{\mathbb{D}}, \mu \in \mathbb{T}$.  Then every
derivation of $\mathcal{T}^{+}(\mathcal{C}_n)$ at $ \pi_{w, \lambda,
\mu}$ is inner if and only if $|\lambda | = 1$. \end{prop}

\begin{proof} Notice that the only primitive cycles in
$\mathcal{C}_n$ are given by \[ e_je_{j+1}\cdots e_ne_1e_2 \cdots
e_{j-1}, \] where $1 \leq j \leq n$.  In this case notice that
$\iota_{w}$ is a cyclic automorphism, in the sense of \cite[Section
2]{Alaimia:1999}. It follows that $\pi_{w, \lambda, \mu}$ factors
through $\mathcal{T}^{+}(\mathcal{C}_n)$ as a completely contractive
automorphism followed by evaluation at $ \lambda$.  Let
$\pi_{\lambda}$ denote the representation of
$\mathcal{T}^{+}(\mathcal{C}_n)$ given by evaluation at $ \lambda$.
It was shown in \cite{Duncan:2005} that every continuous derivation
of $\mathcal{T}^{+}(\mathcal{C}_n)$ at $\pi_{\lambda}$ is inner if
and only if $|\lambda| = 1$. The result now follows.\end{proof}

Notice that there exist nonzero inner derivations for
$\mathcal{T}^{+}(\mathcal{C}_n)$ at $ \pi_{w, \lambda, \mu}$ for all
$ \lambda$ and $ \mu$.  Now, as we did in the case of
representations, we will use the derivations of
$\mathcal{T}^{+}(\mathcal{C}_n)$ to tell us about the derivations of
$\mathcal{T}^{+}(Q)$ for an arbitrary graph.

\begin{dfn} Let $\pi: A \rightarrow B$ be a completely contractive
representation which factors through $C$, via $\tilde{\pi} \circ
\iota$.  We say that a continuous derivation, $D: A \rightarrow B$,
at $\pi$ {\em factors (continuously) through $C$} if there exists a
(continuous) derivation, $\tilde{D}: C \rightarrow B$, at
$\tilde{\pi}$ such that $ D(a) = \tilde{D} \circ \iota (a)$ for all
$ a \in A$. \end{dfn}

If $w$ is a primitive cycle of length $n$ in $Q$, $ \lambda \in
\overline{\mathbb{D}}$, and $\mu \in \mathbb{T}$, recall that the
representation $\pi_{w, \lambda, \mu}:\mathcal{T}^{+}(Q) \rightarrow
M_n$ factors through $\mathcal{T}^{+}(\mathcal{C}_n)$ via the map $
\pi_{\lambda, \mu} \circ \iota_w$. Now, if $D:
\mathcal{T}^{+}(\mathcal{C}_n) \rightarrow M_n$ is a continuous
derivation at $ \pi_{\lambda, \mu}$ then the map $D \circ \iota_w$
induces a derivation on $\mathcal{T}^{+}(Q)$. It is clear that the
induced derivation will factor through
$\mathcal{T}^{+}(\mathcal{C}_n)$.  It follows that the derivations
of $\mathcal{T}^{+}(Q)$ at $ \pi_{w, \lambda, \mu}$ that factor
through $\mathcal{T}^{+}(\mathcal{C}_n)$ are completely determined
by the derivations of $\mathcal{T}^{+}(\mathcal{C}_n)$ at
$\pi_{\lambda, \mu}$. Thus the description of derivations that
factor continuously through $\mathcal{T}^{+}(\mathcal{C}_n)$ can be
easily understood from \cite{Duncan:2005}, where the continuous
point derivations of $ \mathcal{T}^{+}(\mathcal{C}_n)$ are studied.
The following is just a restatement of the description of continuous
derivations of $\mathcal{T}^{+}(\mathcal{C}_n)$ at the
representation given by evaluation of
$\mathcal{T}^{+}(\mathcal{C}_n)$ at $ \lambda$.

\begin{prop} Let $w$ be a primitive cycle of length $n$ in $Q$,
$\lambda \in \overline{\mathbb{D}}$, and $ \mu \in \mathbb{T}$.
Assume that $D: \mathcal{T}^{+}(Q) \rightarrow M_n$ is a non-inner
derivation at $ \pi_{w, \lambda, \mu}$ which factors continuously
through $\mathcal{T}^{+}(\mathcal{C}_n)$, then $| \lambda| < 1$.
\end{prop}

We now determine those derivations that do not factor through
$\mathcal{T}^{+}(\mathcal{C}_n)$. We begin by describing a method to
check whether a derivation factors through
$\mathcal{T}^{+}(\mathcal{C}_n)$.

\begin{prop} Let $Q$ be a directed graph, $w$ a primitive cycle of
length $n$ in $Q$, $ \lambda \in \mathbb{D}$, and $ \mu \in
\mathbb{T}$.  If $D:\mathcal{T}^{+}(Q) \rightarrow M_n$ is a
continuous derivation at $ \pi_{w, \lambda, \mu}$, then $D$ factors
through $\mathcal{T}^{+}(\mathcal{C}_n)$ if and only if $D|_{\ker
\iota_w} \equiv 0$, where $\iota_w$ is the canonical map from
$\mathcal{T}^{+}(Q)$ into
$\mathcal{T}^{+}(\mathcal{C}_n)$.\end{prop}

\begin{proof}  Certainly if $a \in \ker \iota_w$ and $D(a) \neq 0$ then
$D$ cannot factor through $\mathcal{T}^{+}(\mathcal{C}_n)$, else the
induced derivation would send $0$ to a nonzero element of $M_n$.

We next assume that $D|_{\ker \iota_w} = 0$.  Since the range of
$\iota_w$ is complemented as a Banach subspace of
$\mathcal{T}^{+}(\mathcal{C}_n)$, every element $x \in
\mathcal{T}^{+}(\mathcal{C}_n)$ can be written uniquely as $x_r +
x_k$ where $ x_r \in \ran \iota_w$ and $ x_k$ is in the orthogonal
complement of $\ran \iota_w$.  We need only define a continuous
derivation on $\ran \iota_w$ and extend it continuously to
$\mathcal{T}^{+}(\mathcal{C}_n)$ by sending the orthogonal
complement of $\ran \iota_w$ to zero.

This follows by defining a map $\widehat{D}: \ran (\iota_w)
\rightarrow M_n$ via the definition $ \widehat{D}(\iota_w(x)) =
D(x)$. We need only show that this map defines a derivation at $
\tau_{\lambda} \circ \pi_{\mu}$.  Notice that \begin{align*}
\widehat{D}( \iota_w(x) \iota_w(y)) &= \widehat{D}(\iota_w(xy)) \\
&= \pi_{w, \lambda, \mu}(x) D(y) + D(x) \pi_{w, \lambda, \mu}(y)
\\ & = \tau_{\lambda} \circ \pi_{\mu} (\iota_w(x))
\widehat{D}(\iota_w(y)) + \widehat{D}(\iota_w(x)) \tau_{\lambda}
\circ \pi_{\mu} ( \iota_w(y)). \end{align*}  Thus $\widehat{D}$
defines a derivation at $ \tau_{\lambda} \circ \pi_{\mu}$, with $
\widehat{D} \circ \iota_w = D$. \end{proof}

\begin{cor}  For the representation $\pi_{w, \lambda, \mu}: \mathcal{T}^{+}(Q)
\rightarrow M_n$ with $w$ a primitive cycle in $Q$ let D be a
derivation at $ \pi_{w, \lambda, \mu}$.  If $D$ is inner, then $D$
factors continuously through $\mathcal{T}^{+}(\mathcal{C}_n)$, where
$n$ is the length of the primitive cycle $w$.
\end{cor}

\begin{proof}  Since $D$ is inner we know that $D|_{\ker (\pi_{w,
\lambda, \mu})} \equiv 0$.  But notice that $\ker( \iota_w)
\subseteq \ker ( \pi_{w, \lambda, \mu})$ and the result
follows.\end{proof}

Thus, we do not get any new inner derivations on
$\mathcal{T}^{+}(Q)$ besides those obtained via
$\mathcal{T}^{+}(\mathcal{C}_n)$.  We do, however, get non-inner
derivations that do not factor through
$\mathcal{T}^{+}(\mathcal{C}_n)$. We describe some of these
derivations now.

\begin{thm} Let $w = e_1e_2 \cdots e_n$ be a primitive cycle in
$Q$, $\lambda \in \mathbb{D}$, and $ \mu \in \mathbb{T}$. There
exist derivations at $ \pi_{w, \lambda, \mu}: \mathcal{T}^{+}(Q)
\rightarrow M_n$ which do not factor through
$\mathcal{T}^{+}(\mathcal{C}_n)$ if and only if either
\begin{description} \item[(i)] there is an edge $e \neq e_i$
for all $i$ such that there are $j$ and $k$ with $r(e) = s(e_j)$ and
$s(e) = s(e_k).$

\item[ or (ii)]  there is an edge $e_i$ such that $r(e_i) =
s(e_i)$. \end{description}\end{thm}

\begin{proof} Notice that if $ v$ is a vertex with $v \neq s(e_j)$
for all $j$, then \begin{align*} D(P_v) &= D(P_vP_v) \\ &= \pi_{w,
\lambda, \mu}(P_v)D(P_v) + D(P_v) \pi_{w, \lambda, \mu}(P_v) \\ & =
0 D(P_v) + D(P_v) ) \\ & = 0. \end{align*} Now if $e$ is an edge
with $r(e) \neq s(e_j)$ for all $j$, then \begin{align*} D(L_e) &=
D(P_{r(e)}L_e) \\ &= \pi_{w, \lambda, \mu}(P_{r(e)}) D(L_e) +
D(P_{r(e)}) \pi_{w, \lambda, \mu}(L_e) \\ &= 0 D(L_e) + 0 \cdot 0 \\
&= 0. \end{align*}  Similarly, if $ e$ is an edge with $s(e) \neq
s(e_j)$ for all $ j$, then $D(L_e) = 0$. Hence if $D(L_e) \neq 0$
then $s(e) = s(e_j)$ for some $j$ and $ r(e) = s(e_k)$ for some $k$.

Thus, if there are no edges that satisfy case (i) or case (ii), then
for any edge with $ \iota_w (L_e) = 0$, we know that $D(L_e) = 0$.
Similarly if $\iota_w(P_v) = 0$, then $D(P_v) = 0$.  Further if $X$
is in the ideal generated by such $P_v$ and $L_e$ then $D(X) = 0$.
But notice that the ideal generated by such $P_v$ and $L_e$ contains
the ideal $ \ker \iota_w$ and hence $D$ factors through
$\mathcal{T}^{+}(\mathcal{C}_n)$.

For the converse we assume that $ |\lambda |<1$.

Let $e$ be an edge with $e \neq e_i$ for all $i$ such that $
P_{r(e)}$ and $P_{s(e)}$ are not in $ \ker \pi_{w, \lambda, \mu}$.
Notice that the ideal generated by $L_e$ is complemented as a Banach
subspace of $\mathcal{L}_Q$, denote the ideal by $\langle L_e
\rangle$ and the orthogonal complement by $\langle L_e \rangle^c$.
We define a derivation on $\mathcal{L}_Q$ at $ \pi_{w, \lambda,
\mu}$ by first letting $D|_{\langle L_e \rangle^c}$ be constantly
zero. We now define $D(L_e) = e_{j,i}$ where $s(e) = s(e_i)$ and
$r(e) = s(e_j)$. We claim that this induces a continuous derivation
on $\langle L_e \rangle$.  The restriction of this derivation to
$A_Q$ will be a derivation at $\pi_{w, \lambda, \mu}$ which does not
factor through $A_{\mathcal{C}_n}$.  Notice first that $ P_{r(e)}L_e
P_{s(e)} = L_e$ and hence the derivation property implies that
$D(L_e) = e_{j,j}D(L_e) e_{i,i}$ which is satisfied by $e_{j,i}$.

Next notice that every element $X \in \langle L_e \rangle$ can be
written as $X = L_e \tilde{X}$ where $\tilde{X} \in \mathcal{L}_Q$,
see \cite{Jury-Kribs:2005}.  As $L_e$ is an isometry we also know
that $ \| X \| = \|\tilde{X} \|$. By the derivation property we know
that $D(X) = D(L_e) \pi_{w, \lambda, \mu}( \tilde{X})$.  Now
\begin{align*}
\| D(X) \| &= \| D(L_e) \| \| \pi_{w, \lambda, \mu}(\tilde{X}) \| \\
& \leq \| D(L_e) \| \| X \|
\end{align*} and hence $D$ is continuous.  Restricting the derivation
to $\mathcal{T}^{+}(Q)$ gives a continuous non-inner derivation at
$\pi_{w, \lambda, \mu}$.

For the case in which $\pi_{w, \lambda, \mu}(L_e) \neq 0$ but $r(e)
= s(e)$ we use the same proof as in the preceding case except here
we let $D(L_e) = \sum_{r(e_i) = r(e)} e_{i,i}$ to construct the
derivation.  The same proof of continuity will work in this case as
in the previous.  Notice however that $D(L_e^{| \{ i: e_i = e\}|+1
}) \neq 0$ but $ L_e^{| \{ i: e_i = e\}|+1 } \in \ker (\iota_w)$ and
hence the derivation does not factor through
$\mathcal{T}^{+}(\mathcal{C}_n)$.\end{proof}

\begin{rmk} In the previous proof we needed to assume that
$\lambda \in \mathbb{D}$ to use the special structure of ideals of
$\mathcal{L}_Q$ to get continuity of the derivation.  It is possible
that the derivation will also be continuous at $ \pi_{w, \lambda,
\mu}$ with $ | \lambda | = 1$ but we have not been able to construct
a proof of continuity of the described derivation in this case.
\end{rmk}

It might be reasonable to expect a result along the lines of
\cite[Corollary 2]{Duncan:2005}.  However, it is not clear how one
would piece together different copies of
$\mathcal{T}^{+}(\mathcal{C}_n)$ arising from different primitive
cycles to make a general statement about a transitive graph algebra.
The previous proposition also suggests that dealing with loop edges
will complicate the situation.

We close this section by noticing that this approach will provide
little help in dealing with graphs which are not transitive. In fact
there are non transitive graphs which have derivations that are not
inner.  As examples notice that the algebras $\mathcal{A}_{2n}$ of
\cite{Gilfeather:1984}, which have non-inner derivations, can be
viewed, completely isometric isomorphically, as the graph algebra
arising from nontransitive graphs with $2n$ vertices of the form
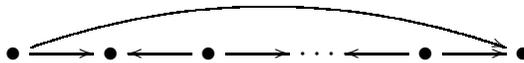
\[ \xymatrix{ \\ {\bullet} \ar[r] \ar@/^1.5pc/[rrrrr] & {\bullet}
& {\bullet} \ar[l] \ar[r] & \cdots & {\bullet} \ar[l] \ar[r] &
{\bullet} } \]

\section{Representations and point derivations into $\mathbb{C}$}

In this section we will deal with point derivations at $\pi$ where
$\pi$ is a multiplicative linear functional on $\mathcal{T}^{+}(Q)$.
In this case, the only inner derivation at $\pi$ is the zero
derivation, since the range of a multiplicative linear functional is
$\mathbb{C}$.  We begin by looking at the multiplicative linear
functionals of $\mathcal{T}^{+}(Q)$ for an arbitrary graph $Q$.
These were described in \cite{Kat-Kribs:2003} as an isomorphism
invariant for the algebra $\mathcal{T}^{+}(Q)$.  The next result is
just a restatement of their description in a manner suitable for our
analysis.

\begin{prop}\label{factoraq} For a directed graph $Q$, let $\pi: \mathcal{T}^{+}(Q) \rightarrow
\mathbb{C}$ be a multiplicative linear functional.  Then there
exists an $n$ with $0 < n \leq \infty$ such that $\pi$ factors
through $A_k$, for all $k \geq n$. \end{prop}

\begin{proof} Notice that as $\pi$ is a representation it will
send projections to projections and hence $\pi(P_v) \in \{ 0, 1 \}$
for all $ v \in V(Q)$.  As the projections $\{ P_v: v \in V(Q) \}$
are orthogonal it follows that if there is a vertex $v_0$ with
$\pi(P_{v_0}) = 1$ then $\pi(P_v) = 0$ for all $ v \in V(Q)
\setminus \{v_0 \}$.  If however, $\pi(P_v) = 0$ for all $ v \in
V(Q)$ then $\pi(L_e) = \pi(L_eP_{s(e)}) = 0$ for all edges $e$ and
hence $\pi$ is identically $0$ and not a multiplicative linear
functional.

Now fix $v_0$ the unique vertex with $\pi(P_{v_0}) = 1$.  Notice
that if $e$ is an edge with either $r(e)$ or $s(e)$ not equal to
$v_0$ then $\pi(L_e) = \pi(P_{r(e)}L_e P_{s(e)}) = 0$.  It follows
that $\pi(L_e) \neq 0$ only if $r(e) = s(e) = v_0$.  Let $n$ be the
number of edges with $r(e) = s(e) = v_0$ and put an label these
edges $e_1, e_2, \cdots , e_n$. Now for $k \geq n$ let $f_1, f_2,
\cdots ,f_k$ be the edges in $B_k$ and define a representation
$\iota: \mathcal{T}^{+}(Q) \rightarrow A_k$ by $ \iota(P_{v}) = 1$
if and only if $ v = v_0$, $\iota(L_{e_i}) = f_i$ for all $i$ and $
\iota(L_e) =0 $ if $ e \neq e_i$.

It is easy to see that this map extends to a completely contractive
representation $\iota: \mathcal{T}^{+}(Q) \rightarrow A_k$.  Now
look at the multiplicative linear functional, call it $\tilde{\pi}:
A_k \rightarrow \mathbb{C}$, which satisfies $\tilde{\pi}(f_i) =
\pi(e_i)$.  This multiplicative linear functional completes the
result.\end{proof}

Next, we notice that the same is true of a derivation at $\pi$,
where $\pi$ is a multiplicative linear functional.

\begin{prop}  For $Q$ a countable directed graph, let $\pi$ be
a multiplicative linear functional for $\mathcal{T}^{+}(Q)$, and
assume that $D: \mathcal{T}^{+}(Q) \rightarrow \mathbb{C}$ is a
continuous derivation at $\pi$. Then there exists an $0 < n \leq
\infty$ such that $D$ factors continuously through $A_n$.
\end{prop}

\begin{proof}  As before there exists a vertex $v_0$ with
$ \pi(P_{v_0}) = 1$ such that $ \{ e: r(e) = s(e) = v_0 \} = \{
e_1,e_2, \cdots, e_n\}$.  recall that $\pi(P_v) = 0$ for any $ v\neq
v_0$.  Now let $\iota: \mathcal{T}^{+}(Q) \rightarrow A_n$ be as in
the preceding proposition.  Notice that $D(P_v) = D(P_vP_v) =
\pi(P_v)D(P_v) + D(P_v)\pi(P_v) = 0$ for all vertices with $v \neq
v_0$.  Similarly if $e$ is an edge with either $r(e) \neq v_0$ or
$s(e) \neq 0$ then $ D(L_e) = D(P_{r(e)}L_e) = D(L_eP_{s(e)})$ one
of which must be zero since $L_e$ and one of $P_{r(e)}$ or
$P_{s(e)}$ must be in $ \ker \pi$.

Now notice that $D(P_{v_0}) = D(P_{v_0}P_{v_0}) = D(P_{v_0})
\pi(P_{v_0}) + \pi(P_{v_0}) D(P_{v_0}) $ and since $\mathbb{C}$ is a
field it follows that $D(P_{v_0}) = 0$.  Define a map $\tilde{D}:
A_n \rightarrow \mathbb{C}$ by first assigning $\tilde{D}(1) = 0$
and $ \tilde{D}(L_{f_i}) = D(L_{e_i})$ and extending using linearity
and the definition of a derivation at $\tilde{\pi}$ (i.e.
$\tilde{D}(xy) = \tilde{D}(x) \tilde{\pi}(y) + \tilde{\pi}(x)
\tilde{D}(y)$).  It is easy to see that $ \tilde{D} \circ \iota (a)
= D(a)$ for all $ a \in \mathcal{T}^{+}(Q)$.  We need only see that
$ \tilde{D}$ is a continuous linear functional.  This, however, is
not difficult as $\| \iota (a) \| = \| a \|$ for all $a$ in the
subalgebra of $\mathcal{T}^{+}(Q)$ generated by $ P_{v_{0}}$ and $
L_{e_i }$. Thus for any $b \in A_n$ there is $ \widehat{b} \in
\mathcal{T}^{+}(Q)$ with $ \iota(\widehat{b}) = b$.  By construction
$ \| \tilde{D}(b) \| =  \| D(\widehat{b}) \| \leq \| D \| \|
\widehat{b} \| = \| D \| \| b \|$ and hence $ \tilde{D}$ is
bounded.\end{proof}

It follows that we need only understand the multiplicative linear
functionals and derivations for the algebras $A_n$.  In what follows
we will use notation as if $n$ is finite since the infinite case
will follow in a manner similar to the finite case.

Recall, from \cite[Theorem 3.3]{Dav-Pitts:1998} or \cite{Pop:1996},
that a multiplicative linear functional $\pi: A_n \rightarrow
\mathbb{C}$ is uniquely given by evaluation at a point $ \lambda = (
\lambda_1, \lambda_2, \cdots, \lambda_n) \in
\overline{\mathbb{B}_n}$ where $\pi(L_{e_i}) = \lambda_i$ for all $1
\leq i \leq n$.  Here, by $\overline{\mathbb{B}_n}$, we mean the
unit ball in $\mathbb{C}_n$ with the usual norm.  We will denote the
multiplicative linear functional arising from evaluation at $
\lambda$ by $ \pi_{\lambda}$.

A further consequence of \cite[Theorem 3.3]{Dav-Pitts:1998} is that
the Gelfand map, $\widehat{a}(\lambda) = \pi_{\lambda}(a)$, induces
a homomorphism from $A_n$ into $A(\mathbb{B}_n)$, the analytic
functions on $\mathbb{B}_n$ with continuous extensions to $
\overline{\mathbb{B}_n}$.

Lastly, notice following the arguments of \cite[Proposition
2.4]{Dav-Pitts:1998} that the commutator ideal of $A_n$, denoted
$\mathfrak{C}_n$, is equal to $\displaystyle{ \bigcap_{ \lambda \in
\overline{\mathbb{B}_n}} \pi_{\lambda}}$.  Putting these facts
together we see that $A_{n, \mathfrak{C}} := A_n / \mathfrak{C}_n$
is a semisimple subalgebra of $A(\mathbb{B}_n)$. This will allow us
to use well known results about point derivations of uniform
algebras to discuss the point derivations of $A_n$.  We summarize
the preceding discussion in the following proposition.

\begin{prop}  Let $ \pi_{\lambda}: A_n \rightarrow \mathbb{C}$ be a
multiplicative linear functional, then $\pi_{\lambda}$ factors
through $A_{n, \mathfrak{C}}$. \end{prop}

A similar result holds for derivations.

\begin{prop} Let $\pi_{\lambda}: A_n \rightarrow \mathbb{C}$ be a
multiplicative linear functional and $D$ be a derivation at $
\pi_{\lambda}$.  Then $D$ factors through $A_{n, \mathfrak{C}}$.
\end{prop}

\begin{proof}  Notice that \begin{align*} D(ab - ba) & = D(ab) -
D(ba) \\ & = (D(a) \pi(b) + \pi(a) D(b)) - (D(b) \pi(a) + \pi(b)
D(a) ) \\ & = 0\end{align*} as $ \mathbb{C}$ is commutative.  It
follows that $ D|_\mathfrak{C} \equiv 0$. Further define $\tilde{D}:
A_{n, \mathfrak{C}} \rightarrow \mathbb{C}$ by $\tilde{D} (a) =
D(\widehat{a})$ where $ \widehat{a} = a + k$ for some $k \in
\mathfrak{C}$.  Notice that $ \tilde{D}$ is well defined since $
D(k) = 0$ for all $ k \in \mathfrak{C}$.

Now there exist $k, k_1, k_2 \in \mathfrak{C}$ with
\begin{align*} \tilde{D}(ab) &= D(\widehat{ab}) \\ & = D(ab+k)
\\ &= D(ab) \\ &= D(a) \pi_{\lambda}(b) + \pi_{\lambda}(a) D(b) \\
&= D(a+k_1) \pi_{\lambda}(b+k_2) + \pi_{\lambda}(a+k_1) D(b+k_2) \\
&= \tilde{D}(a) \tilde{\pi_{\lambda}} \circ \iota (b+k_2) + \tilde{
\pi_{\lambda}} \circ \iota(a+k_1) \tilde{D}(b) \\ &= \tilde{D}(a)
\tilde{\pi_{\lambda}}(b) + \tilde{\pi_{\lambda}}(a) \tilde{D}(b).
\end{align*} Hence $\tilde{D}$ is a derivation at $ \pi_{\lambda}$
satisfying the appropriate property for factoring through $A_n$.
\end{proof}

It is not clear that the factorization must be continuous, we will
see later that it is. On the other hand, we will not need continuity
of the induced derivation for what follows. We now describe the
point derivations of $A_{n, \mathfrak{C}}$ by viewing it as a
subalgebra of $A(\mathbb{B}_n)$. As a corollary we will pull back
the derivations and describe when nontrivial derivations can occur
for $A_n$ at a representation $ \pi_{\lambda}$.

\begin{lem} Let $ \lambda \in \mathbb{B}_n$, then $\pi_{\lambda}:
A_n \rightarrow \mathbb{C}$ factors through $A_{m, \mathfrak{C}}$
where $ m$ is the number of $ \lambda_i$ with $ \lambda_i \neq
0$.\end{lem}

\begin{proof} Let $\{ e_i : 1 \leq i \leq n\}$ be the edges in $B_n$ and
assume without loss of generality that $\pi_{\lambda} (e_i) = 0$ if
and only if $ m+1 \leq i \leq n$.  Now let $\Omega: A_n \rightarrow
A_m$ be the completely contractive representation that sends, for $
1 \leq i \leq m$, $L_{e_i}$ to $L_{f_i}$ where $ \{ f_i \}$ is the
set of edges in $B_m$, and sends $L(e_i)$ to zero when $ m+1 \leq i
\leq n$.  Then notice that $ \pi_{\lambda} $ will factor through
$A_m$ via $\pi_{(\lambda_1, \lambda_2, \cdots , \lambda_m)} \circ
\Omega$.  But now $\pi_{(\lambda_1, \lambda_2, \cdots, \lambda_m)}:
A_m \rightarrow \mathbb{C}$ factors through $A_{m, \mathfrak{C}}$.
Putting the appropriate maps together we get that $ \pi_{\lambda}:
A_n \rightarrow \mathbb{C}$ factors through $A_{m, \mathfrak{C}}$.
\end{proof}

Notice that if $A$ is a commutative operator algebra with identity,
and $\pi:A \rightarrow \mathbb{C}$ is a completely contractive
representation, then $ \ker \pi$ is complemented as a Banach
subspace of $A$.  Denote by $ \overline{( \ker \pi)^2}$ the norm
closure of the ideal generated by $ \{ fg: f, g \in \ker \pi \}$. If
$D: A \rightarrow \mathbb{C}$ is a linear functional such that $D(1)
= 0 $ and $ D|_{\overline{(\ker \pi)^2}} = 0$ then we claim that $D$
is a derivation at $\pi$.  To see this, let $f, g \in A$ and notice
that $(f-\pi(f)\cdot 1)(g- \pi(g)\cdot 1) \in ( \ker \pi)^2$. Hence,
$D((f-\pi(f)\cdot 1)(g- \pi(g)\cdot 1)) = 0$.  Multiplying out and
using linearity of $D$ we get that $D(fg) = D(f) \pi(g) + \pi(f)
D(g)$ and hence $D$ is a derivation at $\pi$.  This argument, which
appears in \cite{Browder:1986}, will be used in the proof of the
next result.

\begin{prop}  If $\pi_{\lambda}(L_i) = 0$ for some $i$, then
there is a unique continuous derivation induced by sending $D(L_i)$
to $1$, and $D(L_k) = 0$ for all $ k \neq i$.\end{prop}

\begin{proof}  Since the range of $\pi_{\lambda}$ is
finite dimensional we know that $\ker(\pi_{\lambda})$ is
complemented as a Banach subspace of $A_n$.  Notice further that
$L_i \in \ker(\pi_{\lambda})$ and yet $ L_i \not\in ( \ker
(\pi_{\lambda}))^2$ and hence the map which sends $L_i$ to $1$ and
all other $L_j$ to zero extends, by Hahn-Banach, to a continuous
derivation on $A_n$. \end{proof}

We say that a derivation at $\pi_{\lambda}$ of this form is the
canonical derivation at $L_i$, denoted $D_i$.  Notice that the
canonical derivation does not factor through $A_{m, \mathfrak{C}}$,
where $m$ is the number of nonzero $\lambda_i$.

\begin{prop} Let $D: A_n \rightarrow \mathbb{C}$ be a continuous
derivation at $\pi_{\lambda}$, and let $m$ be the number of $
\lambda_i$ such that $ \lambda_i \neq 0$, then $D = D_1+ D_2$ where
$D_1$ factors through $A_{m, \mathfrak{C}}$ and $D_2$ is a linear
combination of canonical derivations at $L_j$ where $\lambda_j = 0$.
\end{prop}

\begin{proof} For each $ \lambda_i$ with $ \lambda_i = 0$, let
$\omega_i = D(L_i)$.  Then notice that \[ D_1 = D - \sum_{\lambda_i
= 0} \omega_i \cdot D_i \] is a derivation on $A_n$ such that $D_1$
factors through $A_{m, \mathfrak{C}}$ and the result
follows.\end{proof}

\begin{thm} Let $D: A_n \rightarrow \mathbb{C}$ be
a nontrivial point derivation at $\lambda$ which factors through
$A_{m, \mathfrak{C}}$ where $m$ is the number of $\lambda_i$ with $
\lambda_i \neq 0$.  Then $| \lambda| < 1$. Further the nontrivial
point derivation factors continuously through $A_{m,
\mathfrak{C}}$.\end{thm}

\begin{proof} We will assume without loss of generality that $m
= n$.

Let $f(z_1, z_2, \cdots, z_n)$ denote an arbitrary element of $A_{n,
\mathfrak{C}}$ and assume that $|\lambda| = 1$, where $ \lambda =
(\lambda_1, \lambda_2, \cdots, \lambda_n)$. By assumption $\lambda_1
\neq 0$.  We will see that this implies that any point derivation at
$\lambda$, call it $D$, sends $z_1- \lambda_1$ to zero.  It will
then follow by linearity of $D$ that $D(z_1) = 0.$  A similar
argument will then prove that any point derivation at $\lambda$ is
the zero derivation. Notice that the subalgebra $A:= \{ f \in A_{n,
\mathfrak{C}}: f(z_1, z_2, \cdots, z_n) = f(z_1, \lambda_2, \cdots,
\lambda_n) \}$ is a function algebra in the variable $z_1$.  Further
$z_i - \mu \in A$ and hence $A$ separates the points of $ \{ z_1: \|
(z_1, \lambda_2, \lambda_3, \cdots, \lambda_n)\| \leq 1 \}$.

Now define $\pi: A \rightarrow \mathbb{C}$ by $ \pi(f) =
f(\lambda_1, \lambda_2, \cdots, \lambda_n)$.  This is a
multiplicative linear functional on $A$, and by \cite{Browder:1986}
there does not exist a nonzero point derivation at $\pi$ if there
exists a function $g \in A$ with $\| g(z, \lambda_2, \lambda_3,
\cdots, \lambda_n)\| < \| g(\lambda_1, \lambda_2, \lambda_3, \cdots,
\lambda_n) \|$ for all $ z \neq \lambda_1$.  Notice that the
function $g(z) = \frac{z+\lambda_1}{2 \lambda_1}$ satisfies this
property and hence any point derivation at $\pi$ is the zero
derivation. Now since any point derivation on $A_{m, \mathfrak{C}}$
at $\lambda$ will induce a point derivation on $A$ at $
\pi_{\lambda}$, the point derivation must send $L_1$ to zero.

To see that a nontrivial point derivation factors continuously
through $A_{n, \mathfrak{C}}$ we need only see that any nonzero
point derivation of $A_{n, \mathfrak{C}}$ at $ \lambda$ is unique.
Notice that since $|\lambda | <1$ we know that the representation at
$\lambda$ extends to a $wk^*$ continuous representation of
$\mathcal{L}_n$.  Now notice, from \cite[Theorem 2.10 and Theorem
1.3]{Dav-Pitts:1998}, that the ideal $\ker (\pi_{\lambda}) \subset
\mathcal{L}_n$ is equal to the algebraic ideal generated by $n$
elements, $\{ X_1, X_2, \cdots , X_n \}$.  In particular every
element $a \in A_n$ can be written uniquely as $a_0 + Y_1X_1Z_1 +
Y_2X_2Z_2 + \cdots + Y_nX_nZ_n + a_1$ where $a_0 \not\in \ker
(\pi_{\lambda})$, $ Y_i, X_i \not\in \ker \pi_{\lambda}$ for all
$i$, and $a_1 \in \ker (\pi_{\lambda})^2.$ Now applying the quotient
map $q$ we get a decomposition of every element of $f \in
\mathcal{L}_{n, \mathfrak{C}}$ as $f = f(\lambda) +
\displaystyle{\sum_{i=1}^n q(X_i)g_i(z)} + q(a_1)$ where $g_i(z)
\not \in \ker (\pi_{\lambda})$ and $q(a_1) \in (\ker
(\pi_{\lambda}))^2$.  Notice that any derivation at $ \pi_{\lambda}$
will send $q(a_1)$ to zero.  Now
\begin{align*} D(f) &= D(f(\lambda)) + \sum_{i=1}^n D(q(X_i)g_i(z))
+ D(q(a_1)) \\ &= 0 + \sum_{i=1}^n \left( D(q(X_i)) g_i(\lambda) +
\pi_{\lambda}(X_i)D(g_i(z)) \right) + 0 \\ &= \sum_{i=1}^n
D(q(X_i))g_i(\lambda). \end{align*}  Hence every derivation when
restricted to $A_{n, \mathfrak{C}}$ is a linear combination of
scalar multiples of the continuous derivation that sends $X_i$ to
$1$ and every other $X_j$ to zero.\end{proof}

\begin{cor} Let $D: A_n \rightarrow \mathbb{C}$ be a nontrivial
continuous point derivation at $ \pi_{\lambda}$ then, either
$\lambda_{i} = 0$ for some $i$, or $ |\lambda|< 1$. \end{cor}

Also notice that any derivation that factors through $A_{m,
\mathfrak{C}}$ is unique and the above results describe all point
derivations of $A_n$, and hence of $\mathcal{T}^{+}(Q)$ where $Q$ is
a directed graph. This extends the results of Popescu,
\cite{Pop:1996}, where a description of the point derivations was
given for the representation sending $L_i$ to zero for all $i$.

We close with a result concerning $A_n$-valued derivations of $A_n$.
This is a simple application of an idea in \cite[Theorem 16, p.
92]{Bonsall-Duncan:199x}.

\begin{prop} Let $D: A_n \rightarrow A_n$ be a continuous derivation, then
$D(A_n) \subseteq \mathfrak{C}$. \end{prop}

\begin{proof} For all $z \in \mathbb{C}$ we know that $e^{zD}$ is a
continuous automorphism of $A_n$.  Thus for lambda in $
\overline{\mathbb{B}_n}$ the mapping $ \pi_{\lambda} \circ (e^{zD})$
is a multiplicative linear functional on $A_n$ and hence $ |
\pi_{\lambda}(e^{zD})(a) | \leq \| a \|$ for all $ a \in A_n$.  Now
for $a$ in $A_n$ the mapping $ z \mapsto \pi_{\lambda}(e^{zD})(a) $
is a bounded entire function and hence is constant.  But examining
the power series of this function tells us that the coefficient of
$z$ is $ \pi_{\lambda}(Da) $ which must be zero.  As $ \lambda$ was
arbitrary the result follows. \end{proof}

\begin{rmk} The above proof can be extended to arbitrary graphs
using a characterization of the commutator ideal \cite[Corollary
5.5]{Jury-Kribs:2005} and noting that this ideal is the intersection
of the kernels of all multiplicative linear functionals as is the
case for $A_n$. \end{rmk}

It is of course left open whether one can describe the continuous
$\mathcal{T}^{+}(Q)$-valued derivations of $\mathcal{T}^{+}(Q)$ for
an arbitrary transitive graph $Q$, as was done in \cite{Duncan:2005}
for the graph $\mathcal{C}_n$. We have also not made an attempt to
discuss the higher point cohomology for the graph algebras as an
analogue of the results in in Section 9 of \cite{Johnson:1972}.

\bibliographystyle{plain}

\end{document}